\documentclass[11pt,leqno]{article} 
\usepackage{graphics}

\usepackage{lipsum}

\let\OLDthebibliography\thebibliography
\renewcommand\thebibliography[1]{
  \OLDthebibliography{#1}
  \setlength{\parskip}{2pt}
  \setlength{\itemsep}{2pt plus 0.3ex}
}
\newtheorem{thm}{Theorem}[section]

\newcommand{\beqa}{\begin{eqnarray}}
\newcommand{\eeqa}{\end{eqnarray}}

\newcommand{\beq}{\begin{equation}}
\newcommand{\eeq}{\end{equation}}
\newcommand{\lbl}{\label}
\newcommand{\s}{\; \;}

\newcommand{\la}{\lambda}

\newcommand{\p}{\varphi}

\date{}

\title{Generalized Pohozhaev's identity for radial solutions of $p$-Laplace equations}

\author{
Philip Korman   \\ 
Department of Mathematical Sciences \\ 
University of Cincinnati \\ 
Cincinnati Ohio 45221-0025 \\
\\
}

\begin{document}

\maketitle
\begin{abstract} 
We derive a generalized Pohozhaev's identity for radial solutions of $p$-Laplace equations, by using the approach in \cite{K2},  thus extending the work of H. Br\'{e}zis  and  L. Nirenberg \cite{BN}, where this identity was implicitly used for the Laplace equation.
 \end{abstract}

\begin{flushleft}
\end{flushleft}

\begin{flushleft}
\end{flushleft}

\section{Introduction}
\setcounter{equation}{0}
\setcounter{thm}{0}
\setcounter{lma}{0}

Any solution $u(x)$ of semilinear Dirichlet problem on a bounded smooth domain $D \subset R^n$
\begin{equation}
\label{po1}
\Delta u +f(u)=0 \s \mbox{in $D$}, \s u=0 \s \mbox{on $ \partial D$}
\end{equation}
satisfies the well-known Pohozhaev's identity
\begin{equation}
\label{po2}
\int _{D} \left[2n F(u)+(2-n)uf(u) \right] \, dx=\int _{\partial D} (x \cdot \nu) | \nabla u |^2 \, dS \,.
\end{equation}
Here $F(u)=\int_0^u f(t) \,dt$, and $\nu$ is the unit normal vector on $\partial D$, pointing outside. A standard proof involves multiplication of the equation (\ref{po1}) by $x \cdot \nabla u$ and repeated integration by parts, see e.g., K. Schmitt \cite{S}. In our book \cite{K2} we observed that a more straightforward derivation is to show first that $z \equiv x \cdot \nabla u$ satisfies
\begin{equation}
\label{po3}
\Delta z +f'(u)z=-2f(u) \s \mbox{in $D$}, \s z=0 \s \mbox{on $ \partial D$} \,,
\end{equation}
and then from the equation (\ref{po1}) multiplied by $z$ subtract the equation (\ref{po3}) multiplied by $u$, followed by integration over $D$. We used a similar approach for non-autonomous elliptic systems of Hamiltonian type in \cite{K2} and  \cite{K3}, including systems with power nonlinearities, obtaining an easy derivation of the critical hyperbola, see \cite{K2} for details.
\medskip

For radial solutions on  balls in $R^n$ there is a more general Pohozhaev's identity. It was used implicitly in the classical paper of  H. Brezis and L. Nirenberg \cite{BN}, but it   was not written down in the general form, as presented next. (As above $F(u)=\int_0^u f(t) \, dt$.)
\begin{thm}
Let $u(r) \in C^2[0,1]$ be a solution of 
\begin{equation}
\label{p1}
u'' +\frac{n-1}{r}u'+f(u)=0 \,, \s \mbox{$0<r<1$}\,, \s u'(0)=u(1)=0 \,,
\end{equation}
and let $\psi(r) \in C^2[0,1]$. Then
\begin{equation}
\label{p2}
\s\s \int_0^1 \left[ 2(\psi r^{n-1})' F(u)+\left(2\psi'r^{n-1}-(\psi r^{n-1})' \right)uf(u)-uu'L[\psi]r^{n-3} \right] \, dr
\end{equation}
\[
=\psi (1){u'}^2(1) \,,
\]
where $L[\psi]=r^2\psi''-(n-1)r\psi'+(n-1)\psi$.
\end{thm}

We shall prove a more general $p$-Laplace version of this result, by using the approach described above, and present an application based on \cite{BN}. Similarly to \cite{K2} and \cite{K3} it appears possible to extend these results in two directions: to allow $f(r,u)$ with  $r$ dependence, and to consider systems.
\medskip

Another generalization of radial Pohozhaev's identity, also stimulated by H. Brezis and L. Nirenberg \cite{BN}, was found by F. Catrina \cite{Cat}.

\section{An application}
\setcounter{equation}{0}
\setcounter{thm}{0}
\setcounter{lma}{0}

The generalized Pohozhaev's identity (\ref{p2})  appears to be too involved to use, except in the following three cases: when  $n=3$, or when $\psi(r)=r$, or in case $\psi(r)=r^{n-1}$.
\smallskip

 In case $n=3$,  assuming that $\psi(r) \in C^3[0,1]$ satisfies $\psi(0)=0$, we have $L[\psi]=r^2\psi''-2r\psi'+2\psi$, $L[\psi](0)=0$, and then
\[
-\int_0^1 uu'L[\psi]r^{n-3} \, dr= \frac12 \int_0^1 u^2 \frac{d}{dr} L[\psi] \, dr = \frac12 \int_0^1 u^2 \psi'''r^2 \, dr \,,
\]
and (\ref{p2}) simplifies to become
\[
\int_0^1 \left[ 2(\psi r^{2})' F(u)+\left(2\psi'r^{2}-(\psi r^{2})'\right)uf(u) +\frac{1}{2}u^2 \psi'''r^2  \right] \, dr =\psi (1){u'}^2(1) \,.
\]

\noindent
{\bf \large Example 1 } $f(u)=\la u+u|u|^{p-1}$, with $p \geq 5$. ($5$ is the critical exponent $\frac{n+2}{n-2}$ for $n=3$).
Then $uf(u)=\la u^2+|u|^{p+1}$, $F(u)=\frac{1}{2} \la u^2+\frac{1}{p+1} |u|^{p+1}$, and the last identity becomes
\[
 \int_0^1 \left[ \frac{p+3}{p+1} \psi'r^2-\frac{2(p-1)}{p+1} \psi r \right]|u|^{p+1}     \, dr+\frac{1}{2} \int_0^1 \left(\psi '''+4 \lambda \psi' \right) u^2r^2\, dr=\psi(1){u'}^2(1) \,.
\]
This formula results in a contradiction (proving non-existence of solutions) provided that
\begin{eqnarray}
& \psi(0)=0, \;\; \psi(1) \geq 0 \\ \nonumber
& \psi '''+4 \lambda \psi'=0 \\ \nonumber
& 2(p-1)\psi r-(p+3) \psi'r^2>0 \,. \nonumber
\end{eqnarray}
The equation in the second  line, and the boundary conditions in line one, are satisfied by $\psi (r) = \sin \sqrt{4 \lambda} r$, with $\lambda \in (0,\frac{\pi^2}{4}]$. The last inequality requires that
\[
 \sin \sqrt{4 \lambda} r>\frac{p+3}{2(p-1)} \sqrt{4 \lambda} r \cos \sqrt{4 \lambda} r \,,
\]
or 
\[
 \sin \theta -\gamma  \theta \cos \theta>0\,,
\]
if we denote $ \gamma=\frac{p+3}{2(p-1)}$, and $\theta = \sqrt{4 \lambda} r$. Observe that $\gamma \in (0,1]$, provided that $p \geq 5$, and $\theta  \in (0,\pi)$ for $\lambda \in (0,\frac{\pi^2}{4})$. Then 
\[
 \sin \theta -\gamma  \theta \cos \theta \geq \gamma \left( \sin \theta -  \theta \cos \theta \right)>0 \,.
\]
Conclusion: for  $p \geq 5$, and $\lambda \in (0,\frac{\pi^2}{4}]$ the problem ($n=3$)
\[
u'' +\frac{2}{r}u'+\lambda u+u|u|^{p-1}=0 \,, \s \mbox{$0<r<1$}\,, \s u'(0)=u(1)=0 
\]
has no non-trivial solutions.
\medskip

\noindent
{\bf Remarks}
\begin{enumerate}
  \item The same conclusion holds for other $f(u)$, e.g., for $f(u)=\la u+u|u|^{p-1}+u|u|^{q-1}$, with $q>p \geq 5$.
  \item In case $p>5$ non-existence of solutions for $\la$ small was proved in the same paper  of H. Brezis and L. Nirenberg \cite{BN}, and  in C. Budd  and  J. Norbury \cite{budd}, see also Proposition 1.1 in \cite{K2}.
\end{enumerate}

In case  $p=5$, this example is a part of the classical result of H. Brezis and L. Nirenberg \cite{BN}, who also proved the existence of solutions for $\la \in (\frac{\la _1}{4},\la _1)$ (observe that $\la _1=\pi^2$ for the unit ball in $R^3$). It is remarkable that their non-existence result is sharp for $p=5$.  Let us recall this theorem of H. Brezis and L. Nirenberg \cite{BN} (an extension to sign-changing solutions was later  given in F.V. Atkinson, H. Brezis and L.A. Peletier  \cite{Atki}).
\begin{thm} (\cite{BN})
The problem
\begin{equation}
\label{p3a}
u'' +\frac{2}{r}u'+\lambda u+u^5=0 \,, \s \mbox{$0<r<1$}\,, \s u'(0)=u(1)=0 
\end{equation}
has a positive solution if and only if $\lambda \in (\frac{\la _1}{4},\la _1)$.
\end{thm}

\begin{figure}
\begin{center}
\scalebox{0.68}{\includegraphics{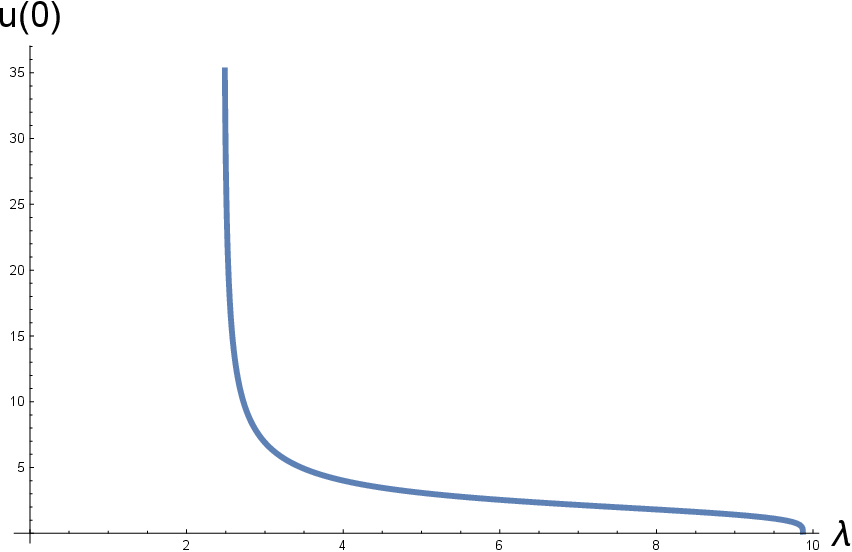}}
\end{center}
\caption{Solution curve of the  Brezis-Nirenberg problem (\ref{p3a})}
\lbl{fig:1}
\end{figure}
\vspace{.1in}

We used {\em Mathematica} to compute the solution curve in the $(\lambda,u(0))$ plane of the  Brezis-Nirenberg problem (\ref{p3a}), presented in Figure \ref{fig:1}, with $u(0)$ giving the maximum value of solutions.  (We used the scaling $u=\sqrt[4]{\lambda} \, z$, to convert this equation into $z'' +\frac{2}{r}z'+\lambda \left( z+z^5 \right)=0$, to which the shoot-and-scale method, described in detail in P. Korman and D.S. Schmidt \cite{KS} applies. A program in {\em Mathematica}  can be downloaded from \cite{KS1}.) The   picture  in Figure \ref{fig:1} indicates that the solution is unique at each $\la$, and in fact  the  uniqueness   follows from the results of M.K. Kwong and   Y. Li  \cite{KL}.
\medskip

Our numerical computations indicate  that the non-existence result on $(0,\frac{\pi^2}{4})$ in Example $1$ is not  sharp for $p>5$, with solution curves tending to infinity at $\la$ larger than $\frac{\pi^2}{4}$. In Figure \ref{fig:2} we present the solution curve of 
\beq
\lbl{bnir2}
u'' +\frac{2}{r}u'+\lambda u+u^6=0 \,, \s \mbox{$0<r<1$}\,, \s u'(0)=u(1)=0 \,.
\eeq
The solution curve has a completely different shape (see \cite{budd} for the asymptotic behavior of this curve), and the smallest value of $\la$ occurs at the first turning point, $\la \approx 5.91> \frac{\pi^2}{4}$.
\medskip

The identity (\ref{p2}) also simplifies in case 
\beq
\lbl{1-1}
L[\psi]=r^2\psi''-(n-1)r\psi'+(n-1)\psi=0 \,.
\eeq
For $n \geq 3$, one solution of this Euler's equation is $\psi=r$, for which   (\ref{p2}) is the classical Pohozhaev's identity:
\[
\int_0^1 \left[ 2nF(u)+(2-n)uf(u) \right] r^{n-1}\, dr={u'}^2(1) \,.
\]
The other solution of (\ref{1-1}) is $\psi=r^{n-1}$, giving
\[
(4n-4) \int_0^1 F(u(r)) r^{2n-3} \, dr={u'}^2(1) \,.
\]
This identity was used by L.A. Peletier  and  J. Serrin \cite{PS2}.
\medskip

In case $n=2$, the solutions of (\ref{1-1}) are $\psi=r$ and $\psi=r \ln r$, leading to similar identities.

\begin{figure}
\begin{center}
\scalebox{0.68}{\includegraphics{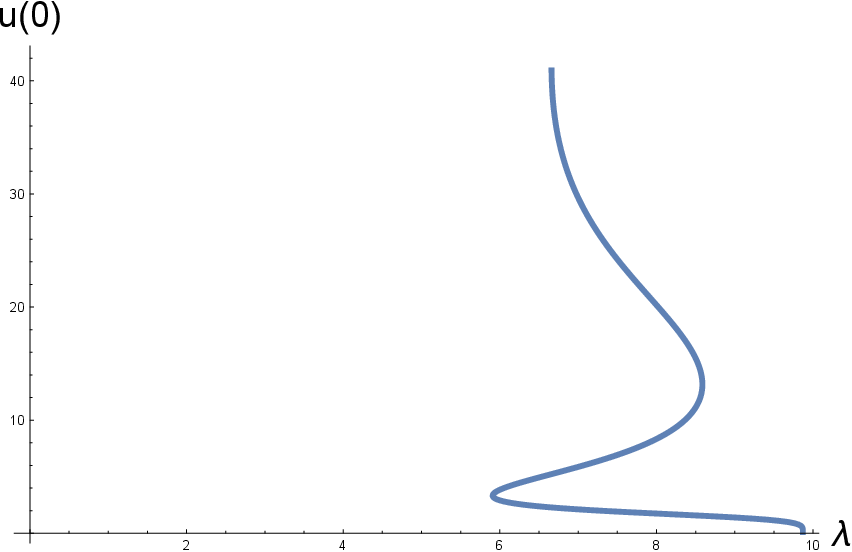}}
\end{center}
\caption{Solution curve of the  supercritical  problem (\ref{bnir2})}
\lbl{fig:2}
\end{figure}
\vspace{.1in}

\section{The $p$-Laplace case}
\setcounter{equation}{0}
\setcounter{thm}{0}
\setcounter{lma}{0}

We present the proof of generalized Pohozhaev's identity next.

\begin{thm}
Let $u(r) \in C^2[0,1]$ be a solution of 
\begin{equation}
\label{1}
\s\s\s \p \left(u'(r)\right)' +\frac{n-1}{r} \p \left(u'(r)\right)+f(u)=0 \s \mbox{$0<r<1$}\,, \s u'(0)=u(1)=0\,,
\end{equation}
with $\p (t)=t|t|^{p-2}$,  $p>1$,
and let $\psi(r) \in C^2[0,1]$. Then
\begin{equation}
\label{2}
\s\s\s \int_0^1 \left[ \left(pF(u)-uf(u) \right) (\psi r^{n-1})' +p  \psi' uf(u)r^{n-1}-\p(u')uL[\psi]r^{n-3} \right] \, dr
\end{equation}
\[
=(p-1) \p(u'(1))\psi (1)u'(1) \,,
\]
where $L[\psi]=(p-1)r^2\psi''-(n-1)r\psi'+(n-1)\psi$.
\end{thm}

\noindent {\bf Proof:} $\s$
Observe that the function $\p (t)$ satisfies
\begin{equation}
\label{3}
t \p' (t)=(p-1) \p (t) \,.
\end{equation}

We claim  that the function $v(r)=\psi (r)u'(r)$ satisfies
\begin{equation}
\label{4}
\left( \p' (u') v'\right)' +\frac{n-1}{r} \p' (u')v'+f'(u)v=-p\psi'f(u)+\frac{\p (u')L[\psi]}{r^2} \,.
\end{equation}
Indeed, using (\ref{3}) and expressing $\p (u')'$ from the equation (\ref{1})
\beqa \nonumber
& \p' (u') v'=\psi ' \p'u'+\psi \p ' (u')u'' =(p-1) \psi ' \p(u') +\psi \p(u')' \\ \nonumber
& =(p-1) \psi ' \p (u') -\frac{n-1}{r} \psi \p (u') -\psi f(u) \,.\\ \nonumber
\eeqa
Then
\beqa \nonumber
& \left( \p' (u') v'\right)'=(p-1) \psi '' \p (u')+(p-1) \psi ' \phi (u')'+\frac{n-1}{r^2} \psi \p (u')-\frac{n-1}{r} \psi' \p (u') \\ \nonumber
& -\frac{n-1}{r} \psi \p (u')'-\psi' f(u)-f'(u)v \,.  \nonumber
\eeqa
Also, using (\ref{3}) again,
\beqa \nonumber
& \frac{n-1}{r} \p '(u') v'=\frac{n-1}{r} \p '(u') \left(\psi 'u'+\psi u'' \right) \\ \nonumber
& =\frac{(n-1)(p-1)}{r}\psi '\p(u') +\frac{n-1}{r} \psi \p (u')' \,. \nonumber
\eeqa
It follows that
\beqa \nonumber
& \left( \p' (u') v'\right)' +\frac{n-1}{r} \p' (u')v'+f'(u)v \\\nonumber
& =(p-1) \psi ''\p +(p-1) \psi ' \left(\p (u')' +\frac{n-1}{r} \p (u') \right) +\frac{n-1}{r^2} \psi \p- \frac{n-1}{r} \psi '\p  \\ \nonumber
& =-p \psi 'f(u) +\p \left[ (p-1)\psi ''- \frac{n-1}{r} \psi '+\frac{n-1}{r^2} \psi\right] \,,\nonumber
\eeqa 
which implies (\ref{4}).
\medskip

Multiplying the equation (\ref{1}) by $(p-1)v$, and subtracting the equation (\ref{4}) multiplied by $u$ gives, in view of (\ref{3}),
\beqa
\label{5}
&  \left[r^{n-1} \left((p-1) \p (u')v- u\p' (u')v' \right) \right]'+r^{n-1}  v \left[(p-1)f(u)-uf'(u) \right] \\
& =pr^{n-1}\psi'uf(u)-r^{n-3}u\p (u')L[\psi] \,. \nonumber
\eeqa
The second term on the left is equal to
\[
\left[pF(u)-uf(u) \right]'\psi r^{n-1}=\left[\left(pF(u)-uf(u) \right) \psi r^{n-1}\right]'-\left(pF(u)-uf(u) \right) \left(\psi r^{n-1} \right)' \,.
\]
Using this identity in (\ref{5}), and integrating over $(0,1)$, we conclude the proof.
\hfill$\diamondsuit$ \medskip

\noindent
{\bf Acknowledgment } It is a pleasure to thank Florin Catrina for a number of useful discussions.

\end{document}